\documentclass[11pt,english]{article}

\usepackage[latin9]{inputenc}
\usepackage{mathtools}
\usepackage{amsmath}
\usepackage{amssymb}
\usepackage[a4paper]{geometry}
\usepackage{setspace}
\makeatletter
\@ifundefined{date}{}{\date{}}
\usepackage{array}
\usepackage{mathrsfs}
\usepackage{color}

\makeatother

\usepackage{babel}
\begin{document}
\title{\textbf{On the Uniqueness of Ein(1) among Linear Combinations of the
Euler-Mascheroni and Euler-Gompertz Constants}}
\author{Michael R. Powers\thanks{Department of Finance, School of Economics and Management, and Schwarzman
College, Tsinghua University, Beijing, China 100084; email: powers@sem.tsinghua.edu.cn.}}
\date{September 30, 2025}
\maketitle
\begin{abstract}
\begin{singlespace}
\noindent From a well-known equation of Hardy, one can derive a simple
linear combination of the Euler-Mascheroni constant ($\gamma=0.577215\ldots$)
and Euler-Gompertz constant ($\delta=0.596347\ldots$): $\gamma+\delta/e=\textrm{Ein}\left(1\right)$.
Although neither $\gamma$ nor $\delta$ is currently known to be
irrational, this linear combination has been shown to be transcendental
(by virtue of the fact that it appears as an algebraic point value
of a particular $E$-function). Moreover, both pairs $(\gamma,\delta)$
and $(\gamma,\delta/e)$ are known to be disjunctively transcendental.
In light of these observations, we investigate the impact of the coefficient
$\alpha$ in combinations of the form $\gamma+\alpha\delta$, and
find that $\alpha=1/e$ is the unique coefficient value such that
canonical Borel-summable divergent series for $\gamma$ and $\delta$
can be linearly combined to force conventional convergence of the
resulting series. We further indicate how this uniqueness property
extends to a sequence of generalized linear combinations, $\gamma^{\left(n\right)}+\alpha\delta^{\left(n\right)}$,
with $\gamma^{\left(n\right)}$ and $\delta^{\left(n\right)}$ given
by (ordinary and conditional) moments of the Gumbel(0,1) probability
distribution.\medskip{}

\end{singlespace}
\end{abstract}
\begin{singlespace}

\section{Introduction}
\end{singlespace}

\begin{singlespace}
\noindent Consider the equation
\begin{equation}
\delta=-e\left[\gamma-{\displaystyle \sum_{k=1}^{\infty}\dfrac{\left(-1\right)^{k+1}}{k\cdot k!}}\right],
\end{equation}
where $\gamma=0.577215\ldots$ denotes the Euler-Mascheroni constant
and $\delta=0.596347\ldots$ denotes the Euler-Gompertz constant.
This equation, introduced by Hardy (1949), can be rewritten as the
following linear combination of $\gamma$ and $\delta$:
\begin{equation}
\gamma+\dfrac{\delta}{e}={\displaystyle \sum_{k=1}^{\infty}\dfrac{\left(-1\right)^{k+1}}{k\cdot k!}}=\textrm{Ein}\left(1\right),
\end{equation}
where $\textrm{Ein}\left(z\right)={\textstyle \sum_{k=1}^{\infty}\left(-1\right)^{k+1}z^{k}/\left(k\cdot k!\right)}$
for $z\in\mathbb{R}$. Both (1) and (2) follow immediately from setting
$z=-1$ in the power-series expansion of the exponential-integral
function,
\begin{equation}
\textrm{Ei}\left(z\right)=\gamma+\ln\left|z\right|+{\displaystyle \sum_{k=1}^{\infty}\dfrac{z^{k}}{k\cdot k!}},\quad z\in\mathbb{R}\setminus\left\{ 0\right\} ,
\end{equation}
and recognizing that $\delta=-e\textrm{Ei}\left(-1\right)$.
\end{singlespace}

\begin{singlespace}
Although neither $\gamma$ nor $\delta$ is known to be irrational,
Rivoal (2012) proved their disjunctive transcendence; that is, at
least one of the two numbers must be transcendental. More recently,
Powers (2025) provided the following probabilistic interpretation
of the elements in (2):
\begin{equation}
\gamma={\displaystyle \int_{-\infty}^{\infty}}x\exp\left(-x-e^{-x}\right)dx=\textrm{E}_{X}\left[X\right],
\end{equation}
\begin{equation}
\delta=-e{\displaystyle \int_{-\infty}^{0}}x\exp\left(-x-e^{-x}\right)dx=-\textrm{E}_{X}\left[X\mid X\leq0\right],
\end{equation}
\begin{equation}
\dfrac{1}{e}={\displaystyle \int_{-\infty}^{0}}\exp\left(-x-e^{-x}\right)dx=\Pr\left\{ X\leq0\right\} ,
\end{equation}
and
\begin{equation}
\sum_{k=1}^{\infty}\dfrac{\left(-1\right)^{k+1}}{k\cdot k!}={\displaystyle \int_{0}^{\infty}}x\exp\left(-x-e^{-x}\right)dx=\textrm{E}_{X}\left[X^{+}\right],
\end{equation}
where $X\sim\textrm{Gumbel}\left(0,1\right)$.\footnote{\noindent The cumulative distribution function of $\textrm{Gumbel}\left(0,1\right)$
is given by $F_{X}\left(x\right)=\exp\left(-e^{-x}\right),\:x\in\mathbb{R}$,
with mean $\gamma$ and variance $\pi^{2}/6$.} He then:

(i) generalized (4), (5), and (7) as
\[
\gamma^{\left(n\right)}={\displaystyle \int_{-\infty}^{\infty}}x^{n}\exp\left(-x-e^{-x}\right)dx=\textrm{E}_{X}\left[X^{n}\right],
\]
\[
\delta^{\left(n\right)}=-e{\displaystyle \int_{-\infty}^{0}}x^{n}\exp\left(-x-e^{-x}\right)dx=-\textrm{E}_{X}\left[X^{n}\mid X\leq0\right],
\]
and
\[
n!{\displaystyle \sum_{k=1}^{\infty}\dfrac{\left(-1\right)^{k+1}}{k^{n}\cdot k!}}={\displaystyle \int_{0}^{\infty}}x^{n}\exp\left(-x-e^{-x}\right)dx=\textrm{E}_{X}\left[\left(X^{+}\right)^{n}\right],
\]
respectively;

(ii) embedded the sequence
\[
\gamma^{\left(n\right)}+\dfrac{\delta^{\left(n\right)}}{e}={\displaystyle \int_{-\infty}^{\infty}}x^{n}\exp\left(-x-e^{-x}\right)dx-{\displaystyle \int_{-\infty}^{0}}x^{n}\exp\left(-x-e^{-x}\right)dx
\]
\begin{equation}
={\displaystyle \int_{0}^{\infty}}x^{n}\exp\left(-x-e^{-x}\right)dx
\end{equation}
for $n\in\mathbb{Z}_{\geq0}$ into the set of $E$-functions
\begin{equation}
F_{n}\left(t\right)={\displaystyle \int_{0}^{\infty}}x^{n}\exp\left(-x-te^{-x}\right)dx;
\end{equation}

(iii) applied the Kolchin-Ostrowski differential-field framework (see
Srinivasan, 2008) to show that, for any $n\in\mathbb{Z}_{\geq1}$,
the set of $F_{0}\left(t\right),F_{1}\left(t\right),\ldots,F_{n}\left(t\right)$
are algebraically independent over $\mathbb{C}\left(t\right)$;

(iv) employed the refined Siegel-Shidlovski theorem (see Beukers,
2006) to prove, by setting $t=1$ in (9), that (8) is transcendental
for all $n\in\mathbb{Z}_{\geq0}$ (including the special case of (2),
for $n=1$); and

(v) concluded further that both pairs $\left(\gamma^{\left(n\right)},\delta^{\left(n\right)}/e\right)$
and $(\gamma^{\left(n\right)},\delta^{\left(n\right)})$ are disjunctively
transcendental for all $n\in\mathbb{Z}_{\geq1}$.

Given that $\gamma+\delta/e=\textrm{Ein}\left(1\right)$ is transcendental,
whereas the arithmetic nature of the linear combination $\gamma+\alpha\delta$
remains unknown for all $\alpha\in\mathbb{R}\setminus\left\{ 1/e\right\} $,
we will investigate the impact of the coefficient $\alpha$ and demonstrate
that the value $\alpha^{*}=1/e$ possesses certain special characteristics.
As shown in Section 2, this coefficient not only equates to the $\textrm{Gumbel}\left(0,1\right)$
probability $\Pr\left\{ X\leq0\right\} $ in (6) (to complete the
probabilistic interpretation of (2)), but also serves as the unique
coefficient for which the linear combination $S_{\gamma}+\alpha S_{\delta}$
-- where $S_{\gamma}$ and $S_{\delta}$ are canonical Borel-summable
divergent series for $\gamma$ and $\delta$, respectively -- converges
conventionally. In Section 3, we discuss the source of this uniqueness
property and indicate how it extends to the sequence of linear combinations,
$\gamma^{\left(n\right)}+\alpha\delta^{\left(n\right)}$, for $n\geq2$.\medskip{}
\medskip{}

\end{singlespace}
\begin{singlespace}

\section{The Uniqueness of $\boldsymbol{\gamma+\delta/e=\textrm{Ein}\left(1\right)}$}
\end{singlespace}

\begin{singlespace}
\noindent In deriving (2) from (3), the coefficient $\alpha=1/e$
in (2) arises directly from writing the Euler-Gompertz constant as
$\delta=-e\textrm{Ei}\left(-1\right)$, which in turn is derived from
the definition $\delta={\textstyle \int_{0}^{\infty}e^{-t}/\left(t+1\right)dt}$.
Although these calculations offer little intuition regarding the provenance
of the indicated coefficient, the probabilistic interpretation given
by (4) through (7) does provide an elegant motivation for $1/e$ (i.e.,
as the value of $\Pr\left\{ X\leq0\right\} $ necessary to join the
partial moment (7) and conditional moment (5) to form the full moment
(4)). Nevertheless, one still can imagine linearly combining $\gamma$
and $\delta$ to achieve something other than $\textrm{E}_{X}\left[X^{+}\right]=\textrm{Ein}\left(1\right)$,
in which case alternative values of $\alpha$ would come into play.
\end{singlespace}

\begin{singlespace}
To show why the choice $\alpha^{*}=1/e$ ($\Longrightarrow\gamma+\alpha\delta=\textrm{Ein}\left(1\right)$)
is unique among all $\alpha\in\mathbb{R}$, we consider its role in
constructing a linear combination of the following two Borel-summable
divergent series:
\begin{equation}
S_{\gamma}\coloneqq{\displaystyle \sum_{k=1}^{\infty}\dfrac{\left(-1\right)^{k}\left(!k\right)}{k}}\overset{\mathcal{B}}{=}\gamma
\end{equation}
and
\begin{equation}
S_{\delta}\coloneqq{\displaystyle \sum_{k=1}^{\infty}\left(-1\right)^{k-1}\left(k-1\right)!}\overset{\mathcal{B}}{=}\delta,
\end{equation}
where $!k=k!{\textstyle \sum_{\ell=0}^{k}}\left(-1\right)^{\ell}/\ell!$
denotes the $k^{\textrm{th}}$ derangement number, ``$\coloneqq$''
means a formal (defined) equality, and ``$\overset{\mathcal{B}}{=}$''
indicates equality via Borel summation.

The divergent alternating-factorial series in (11) is a well-known
Borel-summable expression, whose regularized sum of $\delta$ was
first noted by Euler (see Lagarias, 2013). Moreover, this series famously
was included in the first letter of Srinivasa Ramanujan to G. H. Hardy
(see Berndt and Rankin, 1995). One straightforward means of deriving
(11) is as follows:
\[
\delta={\displaystyle \int_{0}^{\infty}}\dfrac{e^{-t}}{t+1}dt
\]
\[
={\displaystyle \int_{0}^{\infty}}\left(1-t+t^{2}-t^{3}+\cdots\right)e^{-t}dt
\]
\[
=1-\textrm{E}_{T}\left[T\right]+\textrm{E}_{T}\left[T^{2}\right]-\textrm{E}_{T}\left[T^{3}\right]+\cdots
\]
\[
=S_{\delta},
\]
where $T\sim\textrm{Exponential}\left(1\right)$. The Borel-summable
expression for $\gamma$ in (10) does not appear to be widely known
in the literature, but admits a derivation similar to that of (11):
\[
\gamma=-{\displaystyle \int_{0}^{\infty}}\ln\left(t\right)e^{-t}dt
\]
\[
=-{\displaystyle \int_{0}^{\infty}}\left[\left(t-1\right)-\dfrac{\left(t-1\right)^{2}}{2}+\dfrac{\left(t-1\right)^{3}}{3}-\cdots\right]e^{-t}dt
\]
\[
=-\textrm{E}_{T}\left[T-1\right]+\dfrac{1}{2}\textrm{E}_{T}\left[\left(T-1\right)^{2}\right]-\dfrac{1}{3}\textrm{E}_{T}\left[\left(T-1\right)^{3}\right]+\dfrac{1}{4}\textrm{E}_{T}\left[\left(T-1\right)^{4}\right]-\cdots
\]
\[
=S_{\gamma}.
\]

Although there exist infinitely many additional Borel-summable divergent
series for $\gamma$ and $\delta$, those provided in (10) and (11)
are salient in that each is expressed entirely in integer coefficients
and associated with a Borel transform possessing a unique logarithmic
singularity. In particular, note that:

(i) the series (10) has Borel transform 
\[
\mathcal{B}_{\gamma}\left(u\right)=\sum_{k=1}^{\infty}\dfrac{\left(-u\right)^{k}\left(!k\right)}{k\cdot k!}
\]
\[
\Longrightarrow\mathcal{B}_{\gamma}^{\prime}\left(u\right)=\dfrac{1}{u}\sum_{k=1}^{\infty}\dfrac{\left(-u\right)^{k}\left(!k\right)}{k!}
\]
\[
=\dfrac{1}{u}\left[\textrm{E}_{T}\left[e^{-u\left(T-1\right)}\right]-\dfrac{\left(-u\right)^{0}\left(!0\right)}{0!}\right]
\]
\[
=\dfrac{1}{u}\left(e^{u}\textrm{E}_{T}\left[e^{-uT}\right]-1\right)
\]
\[
=\dfrac{1}{u}\left(\dfrac{e^{u}}{u+1}-1\right)
\]
\[
\Longrightarrow\mathcal{B}_{\gamma}\left(u\right)={\displaystyle \int_{0}^{u}}\dfrac{1}{\tau}\left(\dfrac{e^{\tau}}{\tau+1}-1\right)d\tau\quad\textrm{(since }\mathcal{B}_{\gamma}\left(0\right)=0\textrm{)}
\]
\begin{equation}
=\left(\textrm{Ei}\left(u\right)-\ln\left(u\right)\right)+\dfrac{1}{e}\left(\textrm{Ei}\left(1\right)-\textrm{Ei}\left(u+1\right)\right)-\gamma,
\end{equation}
with only a single logarithmic singularity (from $\textrm{Ei}\left(u+1\right)$)
at $u=-1$ (because the singularities from $\textrm{Ei}\left(u\right)$
and $\ln\left(u\right)$ cancel at $u=0$); and

(ii) the series (11) has Borel transform
\[
\mathcal{B}_{\delta}\left(u\right)=u\sum_{k=1}^{\infty}\dfrac{\left(-u\right)^{k-1}\left(k-1\right)!}{k!}
\]
\[
=-\sum_{k=1}^{\infty}\dfrac{\left(-u\right)^{k}}{k}
\]
\begin{equation}
=\ln\left(u+1\right),
\end{equation}
with only a single logarithmic singularity at $u=-1$.

More formally, we take (10) and (11) to be \emph{canonical} divergent
series according to the following criterion: The series is formed
by: (a) expressing the relevant constant ($\gamma$ or $\delta$)
as a Laplace integral, ${\textstyle \int_{0}^{\infty}g\left(t\right)e^{-t}}dt$;
(b) constructing a Taylor-series expansion of the kernel $g\left(t\right)$
at either its minimum positive zero (if one exists) or the boundary
$t=0$; and (c) integrating the Taylor series term by term.

The following proposition identifies the unique impact of the value
$\alpha^{*}=1/e$ on the linear combination $S_{\gamma}+\alpha S_{\delta}$.
\medskip{}

\end{singlespace}

\begin{singlespace}
\noindent\textbf{Proposition 1:} Let $S\left(\alpha\right)\coloneqq S_{\gamma}+\alpha S_{\delta}$,
where $\alpha\in\mathbb{R}$. Then $S\left(\alpha\right)$ is a conventionally
convergent series if and only if $\alpha=\alpha^{*}=1/e$. \medskip{}
\medskip{}

\noindent\textbf{Proof:} Beginning with the expressions for $S_{\gamma}$
and $S_{\delta}$ in (10) and (11), respectively, we find:
\[
S\left(\alpha\right)\coloneqq\sum_{k=1}^{\infty}\dfrac{\left(-1\right)^{k}\left(!k\right)}{k}+\alpha{\displaystyle \sum_{k=1}^{\infty}\left(-1\right)^{k-1}\left(k-1\right)!}
\]
\[
=\sum_{k=2}^{\infty}\dfrac{\left(-1\right)^{k}}{k}k!\left[\sum_{\ell=0}^{k}\dfrac{\left(-1\right)^{\ell}}{\ell!}\right]-\alpha{\displaystyle \sum_{k=2}^{\infty}\left(-1\right)^{k}\left(k-1\right)!}+\alpha
\]
\[
=\sum_{k=2}^{\infty}\left(-1\right)^{k}\left(k-1\right)!\left[\sum_{\ell=0}^{k}\dfrac{\left(-1\right)^{\ell}}{\ell!}\right]-\alpha{\displaystyle \sum_{k=2}^{\infty}\left(-1\right)^{k}\left(k-1\right)!}+\alpha
\]
\[
=\sum_{k=2}^{\infty}\left(-1\right)^{k}\left(k-1\right)!\left[\dfrac{1}{e}-\sum_{\ell=k+1}^{\infty}\dfrac{\left(-1\right)^{\ell}}{\ell!}-\alpha\right]+\alpha
\]
\[
=\left(\dfrac{1}{e}-\alpha\right)\sum_{k=2}^{\infty}\left(-1\right)^{k}\left(k-1\right)!-\sum_{k=2}^{\infty}\sum_{\ell=k+1}^{\infty}\dfrac{\left(-1\right)^{k+\ell}\left(k-1\right)!}{\ell!}+\alpha
\]
\[
=\left(\dfrac{1}{e}-\alpha\right)\sum_{k=2}^{\infty}\left(-1\right)^{k}\left(k-1\right)!-\sum_{k=2}^{\infty}\sum_{m=1}^{\infty}\dfrac{\left(-1\right)^{m}\left(k-1\right)!}{\left(k+m\right)!}+\alpha.
\]
Since the double series is absolutely convergent, we can interchange
the order of summation, obtaining
\[
S\left(\alpha\right)=\left(\dfrac{1}{e}-\alpha\right)\sum_{k=2}^{\infty}\left(-1\right)^{k}\left(k-1\right)!-\sum_{m=1}^{\infty}\left(-1\right)^{m}\left\{ \dfrac{1}{m}\sum_{k=2}^{\infty}\left[\dfrac{\left(k-1\right)!}{\left(k+m-1\right)!}-\dfrac{k!}{\left(k+m\right)!}\right]\right\} +\alpha
\]
\[
=\left(\dfrac{1}{e}-\alpha\right)\sum_{k=2}^{\infty}\left(-1\right)^{k}\left(k-1\right)!-\sum_{m=1}^{\infty}\dfrac{\left(-1\right)^{m}}{m\cdot\left(m+1\right)!}+\alpha
\]
\[
=\left(\dfrac{1}{e}-\alpha\right)\sum_{k=2}^{\infty}\left(-1\right)^{k}\left(k-1\right)!-\sum_{m=1}^{\infty}\dfrac{\left(-1\right)^{m}}{m\cdot m!}-\sum_{m=1}^{\infty}\dfrac{\left(-1\right)^{m+1}}{\left(m+1\right)!}+\alpha
\]
\[
=\left(\dfrac{1}{e}-\alpha\right)\sum_{k=1}^{\infty}\left(-1\right)^{k}\left(k-1\right)!-\sum_{m=1}^{\infty}\dfrac{\left(-1\right)^{m}}{m\cdot m!}
\]
\[
=\left(\alpha-\dfrac{1}{e}\right)S_{\delta}+\sum_{m=1}^{\infty}\dfrac{\left(-1\right)^{m+1}}{m\cdot m!}.
\]
The desired result then follows immediately from the divergence of
$S_{\delta}$ and convergence of\linebreak{}
${\textstyle \sum_{m=1}^{\infty}\left(-1\right)^{m+1}/\left(m\cdot m!\right)=\textrm{Ein}\left(1\right)}$.
$\blacksquare$\medskip{}
\medskip{}

\end{singlespace}
\begin{singlespace}

\section{Discussion}
\end{singlespace}

\begin{singlespace}
\noindent Although Proposition 1 does not directly address the arithmetic
nature of $\gamma+\alpha\delta$ for $\alpha\in\mathbb{R}$, it does
reveal something interesting and important about the case of $\alpha^{*}=1/e$
(for which $\gamma+\delta/e=\textrm{Ein}\left(1\right)$ is known
to be transcendental). Specifically, we see that this particular linear
combination is unique among $\alpha\in\mathbb{R}$ in that it emerges
as a convergent series from $S_{\gamma}+\alpha S_{\delta}$. More
than a simple numerical curiosity, this result sheds light on the
analytical relationship between the two Borel-transform kernels, $\mathcal{B}_{\gamma}\left(u\right)$
and $\mathcal{B}_{\delta}\left(u\right)$, used to define the Euler-Mascheroni
and Euler-Gompertz constants, respectively (via $\gamma={\textstyle \int_{0}^{\infty}}\mathcal{B}_{\gamma}\left(u\right)e^{-u}du$
and $\delta={\textstyle \int_{0}^{\infty}}\mathcal{B}_{\delta}\left(u\right)e^{-u}du$).
\end{singlespace}

\begin{singlespace}
As shown in (12) and (13), each of the kernels is characterized by
a unique logarithmic singularity at $u=-1$, with $\mathcal{B}_{\gamma}\left(u\right)\sim-\left(1/e\right)\ln\left(u+1\right)$
and $\mathcal{B}_{\delta}\left(u\right)\sim\ln\left(u+1\right)$.
Consequently, the associated Stokes constants (for the logarithmic-coefficient
normalization at the singularity) are $-1/e$ and $1$, respectively,
causing the divergent terms in $S_{\gamma}$ and $S_{\delta}$ to
offset one another when combined as $S_{\gamma}+S_{\delta}/e$.

From the Introduction, we recall that Powers (2025) proved that (8)
is transcendental for all $n\in\mathbb{Z}_{\geq0}$, with (2) as a
special case for $n=1$. Therefore, it is worth noting that just as
$\gamma+\delta/e=\textrm{Ein}\left(1\right)$ is unique among linear
combinations $\gamma+\alpha\delta$, the generalization $\gamma^{\left(n\right)}+\delta^{\left(n\right)}/e$
is unique among linear combinations $\gamma^{\left(n\right)}+\alpha\delta^{\left(n\right)}$
for all $n$. Indeed, one can extend $S_{\gamma}$ and $S_{\delta}$
to $n\geq2$ as
\[
S_{\gamma}^{\left(n\right)}\coloneqq\left(-1\right)^{n}n!{\displaystyle \sum_{k=n}^{\infty}\dfrac{s\left(k,n\right)\left(!k\right)}{k!}}\overset{\mathcal{B}}{=}\gamma^{\left(n\right)}
\]
and
\[
S_{\delta}^{\left(n\right)}\coloneqq\left(-1\right)^{n+1}n!{\displaystyle \sum_{k=n}^{\infty}s\left(k,n\right)}\overset{\mathcal{B}}{=}\delta^{\left(n\right)},
\]
respectively, where the $s\left(k,n\right)$ denote signed Stirling
numbers of the first kind. Then, proceeding with an analysis analogous
to that of Section 2, it is possible to show that $S^{\left(n\right)}\left(\alpha\right)\coloneqq S_{\gamma}^{\left(n\right)}+\alpha S_{\delta}^{\left(n\right)}$
is a conventionally convergent series if and only if $\alpha=\alpha^{*}=1/e$.
This is because the Borel-transform kernels are characterized by unique
singularities at $u=-1$, with $\mathcal{B}_{\gamma}^{\left(n\right)}\left(u\right)\sim\left(-1\right)^{n}\left(1/e\right)\left(\ln\left(u+1\right)\right)^{n}$
and $\mathcal{B}_{\delta}^{\left(n\right)}\left(u\right)\sim\left(-1\right)^{n+1}\left(\ln\left(u+1\right)\right)^{n}$.
Consequently, the relevant Stokes constants are $\left(-1\right)^{n}\left(1/e\right)$
and $\left(-1\right)^{n+1}$, respectively, causing the divergent
terms in $S_{\gamma}^{\left(n\right)}$ and $S_{\delta}^{\left(n\right)}$
to offset one another in the linear combination $S_{\gamma}^{\left(n\right)}+S_{\delta}^{\left(n\right)}/e$.\medskip{}
\medskip{}

\end{singlespace}


\begin{thebibliography}{1}
\begin{singlespace}
\bibitem{key-4}Berndt, B. C. and Rankin, R. A., 1995, \emph{Ramanujan:
Letters and Commentary}, American Mathematical Society.

\bibitem{key-1}Beukers, F., 2006, ``A Refined Version of the Siegel-Shidlovskii
Theorem'', \emph{Annals of Mathematics}, 163, 369-379.

\bibitem{key-1}Hardy, G. H., 1949, \emph{Divergent Series}, Clarendon
Press, Oxford, UK.

\bibitem{key-1}Lagarias, J. C., 2013, ``Euler's Constant: Euler's
Work and Modern Developments'', \emph{Bulletin of the American Mathematical
Society}, 50, 4, 527-628.

\bibitem{key-2}Powers, M. R., 2025, ``On the Algebraic Independence
of a Set of Generalized Constants'', arXiv:2508.12123.

\bibitem{key-2}Rivoal, T., 2012, ``On the Arithmetic Nature of the
Values of the Gamma Function, Euler's Constant, and Gompertz's Constant'',
\emph{Michigan Mathematics Journal}, 61, 2, 239-254.

\bibitem{key-1}Srinivasan, V. R., 2008, ``Extensions by Antiderivatives,
Exponentials of Integrals, and by Iterated Logarithms'', arXiv:0811.3004.
\end{singlespace}

\end{thebibliography}
\end{document}